\documentstyle{709}

\begin{document}

\title{INFINITARY AXIOMATIZABILITY OF SLENDER AND COTORSION-FREE GROUPS.}

\author
{Saharon SHELAH\\
Institute of Mathematics, Hebrew University of Jerusalem,\\
Jerusalem, Israel\\
\and
Oren KOLMAN
\thanks{This research was done during the second
author's visit to the Hebrew University in June 1999. 
MSC 2000 classification: primary: 20A15; secondary: 03C55, 03C75, 20A05, 20K20; 
keywords: infinitary logic, axiomatizability, infinite abelian group, slender, 
cotorsion-free, Baer-Specker group.}\\
okolman@member.ams.org}

\maketitle

\begin{abstract}
The classes of slender and cotorsion-free abelian groups are axiomatizable in the
infinitary logics $L_{\infty\omega_1}$ and $L_{\infty\omega}$ respectively. The
Baer-Specker group $\Bbb{Z}^\omega$ is not $L_{\infty\omega_1}$-equivalent to a slender
group.
\end{abstract}

\section{Introduction}
In 1974, Paul Eklof \cite{4} used infinitary logic to generalize some classical theorems
of infinite abelian group theory. He characterized the strongly $\aleph_1$-free groups as
exactly those abelian groups which are $L_{\infty\omega_1}$-equivalent to free abelian
groups, used his criterion to deduce that the class of free abelian groups is not
$L_{\infty\omega_1}$-definable, and showed that the Baer-Specker group $\Bbb{Z}^\omega$
is not $L_{\infty\omega_1}$-equivalent to a free abelian group, strengthening a theorem
of Baer \cite{1} that $\Bbb{Z}^\omega$ is not free. This paper continues in the tradition
allying infinite abelian group theory with infinitary logic. Its main result is the
following:

\begin{theorem}\label{0.1}
The class {\bf SL} of slender abelian groups is axiomatizable in the infinitary logic
$L_{\infty\omega_1}$.
\end{theorem}

\begin{corollary}\label{0.2}
The Baer-Specker group $\Bbb{Z}^\omega$ is not $L_{\infty\omega_1}$-equivalent to a
slender group.
\end{corollary}

Corollary \ref{0.2} improves further the above-mentioned results of Baer \cite{1} and
Eklof \cite{4}. Theorem \ref{0.1} contrasts strikingly with another of Eklof's
corollaries that the class of free groups is not $L_{\infty\omega_1}$-definable. The
proof of Theorem \ref{0.1} also yields that the class of cotorsion-free abelian groups is
$L_{\infty\omega}$-definable.

It is a tacit assumption of the present paper that all the groups are abelian. Any
undefined group-theoretical concepts can be found in Fuchs \cite{8} and Eklof and Mekler
\cite{7}; for logical concepts, see Barwise \cite{2} and Dickmann \cite{3}. We recall
just the essentials and the notation required to read the paper. The {\it Baer-Specker
group} $\Bbb{Z}^\omega$ is the product of countably many copies of $\Bbb{Z}$ (the
integers). Recall that a torsion-free group $G$ is {\it slender} if whenever
$h:\Bbb{Z}^\omega\rightarrow G$ is a homomorphism, then $h(e_n)=0$ for all but finitely
many $n\in \Bbb{N}$, where $e_n$ is the standard $n$-th unit vector
$(0,\dots,0,1,0,\dots)\in\Bbb{Z}^\omega$. We shall use {\bf SL} to denote the class of
slender groups. Since the proof of Theorem \ref{0.1} uses nothing group-theoretical
beyond Nunke's criterion, stated as Theorem \ref{1.1}, the reader wishing to have a
self-contained proof may take Theorem \ref{1.1}(b) as an equivalent definition of
slenderness. Examples of slender groups are $\Bbb{Z}$, and the free groups, since if
$\{G_i:i\in I\}$ is a family of slender groups, then the direct sum $\oplus_{i\in I}G_i$
is slender. Subgroups of slender groups are slender. Strongly $\aleph_1$-free groups are
slender, and hence so are the $\aleph_2$-free groups.

We use $L$ to denote the vocabulary of groups: $L$ contains a binary function symbol $+$
and a constant {\bf O}. For infinite cardinals $\lambda\leq\kappa$, the infinitary
language (or logic) $L_{\kappa\lambda}$ is defined to be the smallest collection of
formulas containing the atomic $L$-formulas and closed under negation, conjunctions
(disjunctions) of sets $\{\varphi_\alpha:\alpha<\mu<\kappa\}$ of formulas, and
existential (universal) quantifier strings of length less than $\lambda$. The language
$L_{\infty\lambda}$ is the union of the languages $L_{\kappa\lambda}$ as $\kappa$ ranges
over the class of infinite cardinals. We write $A\equiv_{\kappa\lambda}B$ for groups $A$
and $B$, to mean that for every $L_{\kappa\lambda}$-sentence $\varphi,\ \varphi$ holds in
$A$ iff $\varphi$ holds in $B$, and in this case, we say that $A$ and $B$ are
$L_{\kappa\lambda}$-equivalent. There is a well-known algebraic characterization of
$L_{\kappa\lambda}$-equivalence using systems of partial isomorphisms (see the references
\cite{2} or \cite{3} for details). A class {\bf C} of groups is
$L_{\infty\lambda}$-definable (or $L_{\infty\lambda}$-axiomatizable) if for some
$L_{\infty\lambda}$-sentence $\varphi$, for every group $G,\ G\in{\bf C}$ iff $\varphi$
holds in $G$.

Many applications of infinitary logic to the study of infinite abelian groups are
described in sections 4, 5 and 7 of the survey \cite{5}; the paper 
\cite{6} contains
further results and bibliographical references.

Group-theoretical interest in $L_{\infty\omega_1}$-axiomatizability stems also in part
from a general thesis of classification theory that measures the complexity of a complete
system of invariants for isomorphism for a class of abelian groups in terms of whether
these invariants can be described by $L_{\infty\omega_1}$-sentences (\cite{12}, p.292;
\cite{11}). By Corollary 1.9 of \cite{4}, there is a group $G(\aleph_1)$ of cardinality
$\aleph_1$ which is $L_{\infty\omega_1}$-equivalent, but not isomorphic, to the free
group $F(\aleph_1)$ of cardinality $\aleph_1$. Since $F(\aleph_1)$ is slender, it follows
that $G(\aleph_1)$ is also slender. In the framework of classification theory, the
results of this paper establish that the first-order theory of the class {\bf SL} does
not have a good structure theorem, since the possible invariants distinguishing the
slender groups $F(\aleph_1)$ and $G(\aleph_1)$ cannot have a simple definition. For the
same reason, the first-order theory of the cotorsion-free groups has no good structure
theorem either.

Finally, let us explain why Theorem \ref{0.1} is the best possible, i.e. why the language
$L_{\infty\omega_1}$ is the right one for axiomatizing {\bf SL}. Note first that {\bf SL}
is not first-order definable, since if $D$ is any non-principal ultrafilter over
$\omega$, then the ultrapower $\Bbb{Z}^\omega/D$ is not slender, but, by $\L$o\v s Lemma,
$\Bbb{Z}^\omega/D\equiv_{\omega\omega}\Bbb{Z}$, and $\Bbb{Z}$ is slender. So {\bf SL} is
not $L_{\omega\omega}$-definable. Nor is {\bf SL} $L_{\infty\omega}$-definable: Barwise
\cite{2} shows that $\Bbb{Z}^\omega\equiv_{\infty\omega}\oplus_{n\in\omega}\Bbb{Z}$, and
again $\oplus_{n\in\omega}\Bbb{Z}$, being free, is slender. On the other hand, it is easy
to see that for every cardinal $\kappa$ greater than the continuum, the class {\bf SL} is
$L_{\infty\kappa}$-definable. So the question of real interest is whether {\bf SL} is
$L_{\infty\kappa}$-definable for $\aleph_1\leq\kappa\leq2^{\aleph_0}$. Theorem \ref{0.1}
and Corollary \ref{0.2} answer this question in the strongest possible form. All our
results are theorems of ZFC (ordinary set theory), and this emphasizes the profound
difference, as far as axiomatizability goes, between the classes of slender (and
cotorsion-free) groups on the one hand, and the free groups on the other.

\section{The proof of Theorem \ref{0.1}}

The key algebraic component in the proof of Theorem \ref{0.1} is the following well-known
characterization of the slender groups, due to Nunke \cite{10} (or see \cite{7} or
\cite{8}).

\begin{theorem}\label{1.1}

For an abelian group $G$, the following conditions are equivalent:
\begin{itemize}
\item[(a)] $G$ is slender;
\item[(b)] $G$ does not contain a subgroup isomorphic to any of the following groups: $\Bbb{Q}$
(the rationals under addition), $\Bbb{Z}^\omega,\ \Bbb{Z}(p)$ (the cyclic group of order
$p$) or $J_p$ (the $p$-adic integers) for any prime $p$.
\end{itemize}
\end{theorem}

It will be convenient to split the proof of Theorem \ref{0.1} into some facts and
propositions which will all have the general form: for each group $B$ in Theorem
\ref{1.1}(b) there is a sentence $\varphi$ such that for every torsion-free group $G,\
\varphi$ is true in $G$ iff $G$ has a subgroup isomorphic to $B$.

Remembering that $\Bbb{Q}$ has a generating system $\{q_n:n\in\omega\}$ with
$q_n=(n+1)q_{n+1}$, one verifies easily the following fact.

\begin{fact}\label{1.2} For any torsion-free group $G,\ G$ contains an isomorphic copy of
$\Bbb{Q}$ iff $G$ satisfies the $L_{\omega_1\omega}$-sentence $\psi^{\Bbb
Q}\equiv(\exists x)\big(x\ne0\wedge\bigwedge_{n\in\omega}(n!$ divides $x$)\big).
\end{fact}

Next we deal with the problem of expressing the assertion "$G$ contains an isomorphic
copy of $\Bbb{Z}^\omega$" by an $L_{\infty\omega_1}$-sentence. Let
$S^*=\{a\in\Bbb{Z}^\omega:n!$ divides $a_n\}$. Note that $S^*$ is a subgroup of
$\Bbb{Z}^\omega$. For $a\in S^*$, let $\varphi_a(y,x_0,x_1,\dots,x_n,\dots)$ be the
following $L_{\omega_1\omega_1}$-formula:
\[
\bigwedge_{n\in\omega}\big(n!\ {\rm divides}\ (y-\Sigma_{l<n}a(l)x_l)\big). \
\]

\begin{fact}\label{1.3} Suppose that $G$ is torsion-free and contains no isomorphic copy
of $\Bbb{Q}$. Then:
\begin{equation}\label{1.3.1}
G\models\ \forall\ (x_0,x_1,\dots,x_n,\dots)\exists^{\leq1}y\varphi_a(y,x_0,x_1,\dots,
x_n,\dots);
\end{equation}
if $h$ is an embedding of $\Bbb{Z}^\omega$ into $G$, then for every $a\in S^*$
\begin{equation}\label{1.3.2}
G\models\varphi_a[h(a),h(e_0),\dots,h(e_n),\dots]
\end{equation}
\end{fact}

\begin{proof}
To see that claim (\ref{1.3.1}) is true, suppose otherwise, and let $c_n\ (n\in\omega)$
and $d_1\ne d_2$ in $G$ witness its failure. Then for every $n\in\omega,\ \bigwedge_{i=
1,2}G\models\big(n!$ divides $(d_i-\Sigma_{l<n}a(l)c_l)\big)$. Subtracting, it follows
that $G\models(d_1-d_2\ne0)\wedge\bigwedge_{n\in\omega}\big(n!$ divides $(d_1-d_2)\big)$,
and hence $G$ satisfies $\psi^{\Bbb Q}$. By Fact \ref{1.2}, this contradicts the
hypothesis that $G$ contains no isomorphic copy of $\Bbb{Q}$.

Claim (\ref{1.3.2}) is immediate because $h$ is an embedding, and the elements $a,e_0,
e_1,\dots,e_n,\dots$ satisfy the formula $\varphi_a(y,x_0,x_1,\dots,x_n,\dots)$ in $\Bbb
{Z}^\omega$.

Let $\varphi^*(x_0,x_1,\dots,x_n,\dots)$ be the logical conjunction of the following
formulas:
\begin{itemize}
\item[(A)]$\bigwedge\big\{(\exists y)\varphi_a(y,x_0,x_1,\dots,x_n,\dots):a\in S^*\big\}$
\item[(B)]$\bigwedge\big\{(\exists y_1\exists y_2\exists y_3)\big(\bigwedge_{i=1,2,3}\varphi_{a^i}
(y_i,x_0,x_1,\dots,x_n,\dots)\wedge(y_3=y_1-y_2)\big):a^1,a^2,a^3\in S^*,a^3=a^1-a^2\}$
\item[(C)]$\bigwedge_{m\in\omega}(x_0,x_1,\dots,x_{m-1}$ are independent)
\item[(D)]$\bigwedge\big\{(\exists y)\big(\varphi_a(y,x_0,x_1,\dots,x_n,\dots)\wedge y\ne0\big)
:a\in S^*,a\ne0\big\}$.
\end{itemize}

Observe that $\varphi^*(x_0,x_1,\dots,x_n,\dots)$ is an $L_{(2^{\aleph_0})^+
\aleph_1}$-formula. Let $\varphi^{\Bbb{Z}^\omega}$ be the $L_{(2^{\aleph_0})^+\aleph_
1}$-sentence $\exists(x_0,x_1,\dots,x_n,\dots)\varphi^*(x_0,x_1,\dots,x_n,\dots)$.
\end{proof}

\begin{proposition}\label{1.4} Suppose that $G$ is torsion-free and contains no subgroup
isomorphic to $\Bbb{Q}$. The following conditions are equivalent:
\begin{itemize}
\item[(\ref{1.4}.1)] there exists an embedding $h:\Bbb{Z}^\omega\rightarrow G$;
\item[(\ref{1.4}.2)] $G\models\varphi^{\Bbb{Z}^\omega}$.
\end{itemize}
\end{proposition}

\begin{proof}
The implication (\ref{1.4}.1) $\Rightarrow$ (\ref{1.4}.2) is easy, since $\Bbb{Z}^\omega
\models\varphi^*[e_0,e_1,\dots,e_n,\dots]$, and one can use the embedding $h$ to find
witnesses in $G$ for the existential quantifiers of the conjuncts of $\varphi^*$.

Conversely, suppose that (\ref{1.4}.2) holds, and $c_0,c_1,\dots,c_n,\dots$ are elements
of $G$ satisfying $G\models\varphi^*[c_0,c_1,\dots,c_n,\dots]$. The map $h_0:\Bbb{Z}^
\omega\rightarrow S^*$ given by $h_0(a)=\langle\dots,n!a_n,\dots\rangle$ is an embedding,
so to establish (\ref{1.4}.1), it is enough to find an embedding $d:S^*\rightarrow G$
(and then $h=d\circ h_0$ embeds $\Bbb{Z}^\omega$ into $G$). For an element $a\in S^*$,
since $G\models\varphi^*[c_0,c_1,\dots,c_n,\dots]$, there exists $d_a\in G$ such that $G
\models\varphi_a[d_a,c_0,c_1,\dots,c_n,\dots]$. By Fact \ref{1.3.1}, the element $d_a$ is
uniquely defined. Claim: the map $d:S^*\rightarrow G$ given by $d(a)=d_a$ is an
embedding. This is easy to show: to see that $d$ is a homomorphism, note that e.g., for
$a,b\in S^*$, for every $n\in\omega,\ G\models\big(n!$ divides $\big(d_{(a+b)}-(d_a+
d_b)\big)\big)$, and hence $d(a+b)=d(a)+d(b)$ (otherwise Fact \ref{1.2} will reveal a
copy of $\Bbb{Q}$ in $G$). A similar argument using the appropriate conjunct of (D) in
$\varphi^*$ proves that $d$ is one-to-one.
\end{proof}

The cases $\Bbb{Z}(p)$ are trivial. Since the group table of $\Bbb{Z}(p)$ can be
completely described by an $L_{\omega\omega}$-sentence, for every prime $p$, there is an
$L_{\omega\omega}$-sentence $\varphi^{\Bbb{Z}(p)}$ such that:

\begin{fact}\label{1.5} For every group $G,\ \varphi^{\Bbb{Z}(p)}$ is true in $G$ iff $G$
has a subgroup isomorphic to $\Bbb{Z}(p)$.
\end{fact}

To deal with the p-adic integers $J_p$, we shall use some elementary properties
summarized below (see [F1]).
\begin{itemize}
\item[(1)] every $t\in J_p$ can be represented as a formal sum $\Sigma_{j<\omega}s_j\cdot p^j$,
where the coefficients $s_j$ belong to $\{0,1,\dots,p-1\}$.

Since $\{0,1,\dots,p-1\}$ can be replaced by any complete set $\{t_0,t_1,\dots,t_{p-1}\}$
of representatives of $\Bbb{Q}_p$ mod $p\Bbb{Q}_p$ (where $\Bbb{Q}_p$ is the ring of
rational numbers whose denominators are prime to $p$), it follows that:

\item[(2)] if $q$ is a prime different from $p$, then for every $t\in J_p,\ q$ divides $t$ in $J_p$.
\end{itemize}

Also, $J_p$ has cardinality continuum.

We define, for each $t=\Sigma_{j<\omega}s_j\cdot p_j\in J_p$, the
$L_{\omega_1\omega}$-formula $\psi_{p,\ t}(y,x)$ to be the logical conjunction of the
following:
\begin{itemize}
\item[(E)] $\bigwedge_{n\in\omega}\big(p^n$ divides $(y-\Sigma_{j<n}s_j\cdot p_j\cdot x)\big)$
\item[(F)] $\bigwedge\{q$ divides $y$ and $q$ divides $x:q\in\Bbb{Z}, (p,q)=1\}$.
\end{itemize}

Since $J_p\models\psi_{p,\ t}[t,1]$, one obtains immediately:

\begin{fact}\label{1.6} If $h:J_p\rightarrow G$ is an embedding, then for every $t\in
J_p,\ G\models\psi_{p,\ t}[h(t),h(1)]$.
\end{fact}

Let $\psi_p(x)$ be the logical conjunction of the formulas:
\begin{itemize}
\item[(G)] $\bigwedge\{(\exists!y)\psi_{p,\ t}(y,x):t\in J_p\}$
\item[(H)] $\bigwedge\{(\exists y_1\exists y_2\exists y_3)\big(\bigwedge_{i=1,2,3}\psi_{p,\ t^i}
(y_i,x)\wedge(y_3=y_1-y_2)\big):t^1,t^2,t^3\in J_p,\ t^3=t_1-t_2\}$
\item[(K)] $\bigwedge\{\exists y)(\psi_{p,\ t}(y,x)\wedge y\ne0):t\in J_p,t\ne0\}$.
\end{itemize}

Note that $\psi_p(x)$ is an $L_{(2^{\aleph_0})^+\aleph_0}$-formula.

\begin{proposition}\label{1.7} Suppose $G$ is torsion-free and does not contain an
isomorphic copy of $\Bbb{Q}$; let $g\in G\backslash\{0\}$. The following conditions are
equivalent:
\begin{itemize}
\item[(\ref{1.7}.1)] there is an embedding $h:J_p\rightarrow G$ such that $h(1)=g$;
\item[(\ref{1.7}.2)] $G\models\psi_p[g]$.
\end{itemize}
\end{proposition}

\begin{proof}
The implication (\ref{1.7}.2) $\Rightarrow$ (\ref{1.7}.1) is analogous to (\ref{1.4}.2)
$\Rightarrow$ (\ref{1.4}.1).

As regards the direction (\ref{1.7}.1) $\Rightarrow$ (\ref{1.7}.2), we sketch only the
slightly different points, arising from $J_p$, concerning the formula $(\exists!y)
\psi_{p,\ t}(y,x)$. By Fact \ref{1.6}, $G\models\psi_{p,\ t}[h(t),h(1)]$. If for some
$d_1\ne d_2\in G,\ \bigwedge_{i=1,2}G\models\psi_{p,\ t}[d_i,h(1)]$, then, for every
$n\in \omega,\ G\models(p^n$ divides $d_1-d_2)$, since $\bigwedge_{i=1,2}G\models\big(
p^n$ divides $(d_i-\Sigma_{j<n}s_j\cdot p^j\cdot h(1)\big)$. But also, if $q\in\Bbb{Z},
(p,q)=1$, then $\bigwedge_{i=1,2}G\models(q$ divides $d_i)$. Hence $G$ satisfies
$\psi^{\Bbb Q}$, contradicting the hypothesis that $G$ does not contain a copy of $\Bbb
Q$. Thus $h(t)\in G$ is unique such that $G\models\psi_{p,\ t}[h(t),g]$. The remaining
conjuncts of (H) and (K) hold in $G$ because $h$ is an embedding. So $G\models\psi_p[g]$.
\end{proof}

Let $\psi^{J_p}$ be the $L_{(2^{\aleph_0})^+\aleph_0}$-sentence $(\exists x)\psi_p(x)$.

\begin{theorem}\label{1.8} There is an $L_{(2^{\aleph_0})^+\aleph_1}$-sentence $\Phi$ such
that for any group $G$, the following conditions are equivalent:
\begin{itemize}
\item[(\ref{1.8}.1)] $G$ is slender;
\item[(\ref{1.8}.2)] $G\models\Phi$.
\end{itemize}
\end{theorem}

\begin{proof} By Theorem \ref{1.1}, together with Facts \ref{1.2} and \ref{1.5}, and
Propositions \ref{1.4} and \ref{1.7}, the following sentence $\Phi$ is as required:
\[
\neg\big(\psi^{\Bbb Q}\vee(\varphi^{\Bbb{Z}^\omega})\vee(\bigvee_{p\ {\rm prime}}\varphi^{\Bbb{Z}(p)})\vee(\bigvee_{p\
{\rm prime}}\psi^{J_p})\big).
\]

Theorem \ref{0.1} and Corollary \ref{0.2} follow immediately from Theorem \ref{1.8}.
\end{proof}

For the definitions and basic properties of Whitehead and Shelah groups, see \cite{7}.

\begin{corollary}\label{1.9} If $G$ is a $W$-group, then $\Bbb{Z}^\omega$ is not
$L_{(2^{\aleph_0})^+\aleph_1}$-equivalent to $G$.
\end{corollary}

\begin{proof} $G$ is slender \cite{7}.
\end{proof}

Since every Shelah group is slender, one also obtains:

\begin{corollary}\label{1.10} If $G$ is a Shelah group, then $\Bbb{Z}^\omega$ is not
$L_{(2^{\aleph_0})^+\aleph_1}$-equivalent to $G$.
\end{corollary}

>From Theorem \ref{0.1} it follows too that $\Bbb{Z}^\omega$ is not $L_{(2^{\aleph_0})^+
\aleph_1}$-equivalent to an $\aleph_2$-free group (this is already a consequence of Eklof
\cite{4}).

\section{Cotorsion-free groups}

Recall that a group $A$ is {\it cotorsion} if Ext$(J,A)=0$ for all torsion-free groups
$J$. A group $G$ is {\it cotorsion-free} if $G$ does not contain any non-zero subgroup 
which is cotorsion. The following theorem characterizes the cotorsion-free groups (\cite{9};
see \cite{7}):

\begin{theorem}\label{2.1} For any group $G$, the following are equivalent:
\begin{itemize}
\item[(\ref{2.1}.1)] $G$ is cotorsion-free;
\item[(\ref{2.1}.2)] $G$ does not contain a copy of $\Bbb{Q}, \Bbb{Z}(p)$, or $J_p$ for any
prime $p$.
\end{itemize}
\end{theorem}

The sentence $\neg\big(\psi^{\Bbb Q}\vee(\bigvee_{p\ {\rm prime}}\varphi^{\Bbb{Z}(p)})
\vee(\bigvee_{p\ {\rm prime}}\psi^{J_p})\big)$ belongs to the logic $L_{(2^{\aleph_0})^+
\aleph_0}$, and so the following corollary is evident:

\begin{corollary}\label{2.2} The class of cotorsion-free groups is $L_{\infty\omega}$-definable.
\end{corollary}

Since the $L_{\infty\omega_1}$-equivalent, non-isomorphic groups $F(\aleph_1)$ and
$G(\aleph_1)$ of cardinality $\aleph_1$ mentioned in the introductory section of this
paper are cotorsion-free, the first-order theory of the cotorsion-free groups has no good
structure theorem either.

\end{document}